\documentclass[12pt]{article}

\usepackage{amsmath, epsfig, cite}
\usepackage{amssymb}
\usepackage{amsfonts}
\usepackage{latexsym}
\usepackage{float}

\newtheorem{thm}{Theorem}[section]

\newtheorem{conj}[thm]{Conjecture}
\newtheorem{lem}[thm]{Lemma}

\numberwithin{equation}{section}

\newcommand*{\pFq}[5]{{}_{#1}F_{#2}\left[ \begin{matrix} #3\\[5pt] #4\end{matrix};#5\right]}

\begin{document}

\begin{center}
{\Large\bf Supercongruences arising from hypergeometric\\[5pt]
series identities
}
\end{center}

\vskip 2mm \centerline{Ji-Cai Liu}
\begin{center}
{\footnotesize Department of Mathematics, Wenzhou University, Wenzhou 325035, PR China\\
{\tt jcliu2016@gmail.com} }
\end{center}


\vskip 0.7cm \noindent{\bf Abstract.} By using some hypergeometric series identities, we prove two supercongruences on truncated hypergeometric series, one of which is related to a modular Calabi--Yau
threefold, and the other is regarded as $p$-adic analogue of an identity due to Ramanujan.

\vskip 2mm
\noindent{\it MR Subject Classifications}: 11A07, 05A19, 33C20
\section{Introduction}
Let
\begin{align*}
f(z):=\eta^4(2z)\eta^4(4z)=\sum_{n=1}^{\infty}a(n)q^n,
\end{align*}
where $q=e^{2\pi iz}$ and the Dedekind eta function
is given by
\begin{align*}
\eta(z)=q^{\frac{1}{24}}\prod_{n=1}^{\infty}(1-q^n).
\end{align*}
For odd primes $p$, let $N(p)$ denote the number of solutions to the modular
Calabi--Yau threefold:
\begin{align*}
x+\frac{1}{x}+y+\frac{1}{y}+z+\frac{1}{z}+w+\frac{1}{w}=0
\end{align*}
over the finite field with $p$ elements. Ahlgren and Ono \cite{an-mm-2000}, van Geemen and Nygaard \cite{gn-ibid-1995}, and Verrill \cite{verrill-b-1996} showed by different methods that
\begin{align*}
a(p)=p^3-2p^2-7-N(p).
\end{align*}

In 2006, Kilbourn \cite{kilbourn-aa-2006} proved that for any odd prime $p$,
\begin{align}
a(p)\equiv \pFq{4}{3}{\frac{1}{2},&\frac{1}{2},&\frac{1}{2},&\frac{1}{2}}{&1,&1,&1}{1}_{\frac{p-1}{2}}\pmod{p^3}.\label{a1}
\end{align}
Here the truncated hypergeometric series are given by
\begin{align*}
\pFq{r}{s}{a_1,a_2,\cdots,a_r}{b_1,b_2,\cdots,b_s}{z}_n=\sum_{k=0}^n\frac{(a_1)_k (a_2)_k\cdots (a_r)_k}{(b_1)_k (b_2)_k\cdots (b_s)_k}\cdot \frac{z^k}{k!},
\end{align*}
where $(a)_0=1$ and $(a)_k=a(a+1)\cdots (a+k-1)$ for $k\ge 1$.

The first aim of this paper is to prove another supercongruence for $a(p)$.
\begin{thm}\label{t1}
For any prime $p\ge 5$, we have
\begin{align}
a(p)\equiv p\cdot\pFq{4}{3}{\frac{1}{2},&\frac{1}{2},&\frac{1}{2},&\frac{1}{2}}{&1,&\frac{3}{4},&\frac{5}{4}}{1}_{\frac{p-1}{2}} \pmod{p^3}.\label{a2}
\end{align}
\end{thm}

In 1997, Van Hamme \cite[(A.2)]{vanHamme-1997} proposed the following supercongruence conjecture.
\begin{conj} (Van Hamme, 1997)
For any odd prime $p$, we have
\begin{align}
&\pFq{6}{5}{\frac{5}{4},&\frac{1}{2},&\frac{1}{2},&\frac{1}{2},&\frac{1}{2},&\frac{1}{2}}
{&\frac{1}{4},&1,&1,&1,&1}{-1}_{\frac{p-1}{2}}\notag\\
&\equiv
\begin{cases}
-p\Gamma_p\left(\frac{1}{4}\right)^4\pmod{p^3}\quad &\text{if $p\equiv 1\pmod{4}$,}\\[10pt]
0\pmod{p^3}\quad &\text{if $p\equiv 3\pmod{4}$,}\label{a3}
\end{cases}
\end{align}
where $\Gamma_p(\cdot)$ denotes the $p$-adic Gamma function.
\end{conj}

The above supercongruence was regarded as $p$-adic analogue of the following identity due to Ramanujan (announced in his second letter to Hardy on February 27):
\begin{align*}
\pFq{6}{5}{\frac{5}{4},&\frac{1}{2},&\frac{1}{2},&\frac{1}{2},&\frac{1}{2},&\frac{1}{2}}
{&\frac{1}{4},&1,&1,&1,&1}{-1}
=\frac{2}{\Gamma\left(\frac{3}{4}\right)^4},
\end{align*}
which was later proved by Hardy \cite{hardy-plms-1924} and Watson \cite{watson-jlms-1931}. The supercongruence \eqref{a3} was first confirmed by McCarthy and Osburn \cite{mo-2008}.

In 2015, Swisher \cite[Theorem 1.5]{swisher-2015} also showed that \eqref{a3} holds modulo $p^5$ for primes $p\equiv 1\pmod{4}$. Recently, Guo and Schlosser \cite[Theorem 2.2]{gs-2018} established an interesting $q$-analogue of a supercongruence closely related to \eqref{a3}.
By using the software package {\tt Sigma} due to Schneider \cite{schneider-1999}, the author \cite[Theorem 1.3]{liu-jmaa-2018} extended the case $p\equiv 3\pmod{4}$ in \eqref{a3} as follows.
\begin{thm}\label{t2}
Let $p\ge 5$ be a prime. For $p\equiv 3\pmod{4}$, we have
\begin{align}
\pFq{6}{5}{\frac{5}{4},&\frac{1}{2},&\frac{1}{2},&\frac{1}{2},&\frac{1}{2},&\frac{1}{2}}
{&\frac{1}{4},&1,&1,&1,&1}{-1}_{\frac{p-1}{2}}\equiv
-\frac{p^3}{16}\Gamma_p\left(\frac{1}{4}\right)^4\pmod{p^4}.\label{a4}
\end{align}
\end{thm}

However, the proof of \eqref{a4} in \cite{liu-jmaa-2018} is based on software package and seems unnatural. The second aim of this paper is to provide a human proof of \eqref{a4} by hypergeometric series identities, which seems to be more natural.

The rest of this paper is organized as follows. Section 2 is devoted to recalling some properties of Gamma function and $p$-adic Gamma function. We prove Theorems \ref{t1} and \ref{t2} in Sections 3 and 4, respectively.

\section{Preliminary results}
We first recall some properties of Gamma function. The Gamma function $\Gamma(z)$ is an extension of the factorial function, which satisfies the functional equation:
\begin{align}
\Gamma(z+1)=z\Gamma(z).\label{e1}
\end{align}
From the above equation, we immediately deduce that for complex numbers $z$ and positive integers $n$,
\begin{align}
(z)_n=\frac{\Gamma(z+n)}{\Gamma(z)}.\label{e2}
\end{align}
It also satisfies the following reflection formula and duplication formula:
\begin{align}
&\Gamma(z)\Gamma(1-z)=\frac{\pi}{\sin(\pi z)},\label{e3}\\
&\Gamma(z)\Gamma\left(z+\frac{1}{2}\right)=2^{1-2z}\sqrt{\pi}\Gamma(2z).\label{e4}
\end{align}

We next recall the definition and some basic properties of $p$-adic Gamma function. For more details, we refer to \cite[Section 11.6]{cohen-2007}.
Let $p$ be an odd prime and $\mathbb{Z}_p$ denote the set of all $p$-adic integers. For
$x\in \mathbb{Z}_p$, the $p$-adic Gamma function is defined as
\begin{align*}
\Gamma_p(x)=\lim_{m\to x}(-1)^m\prod_{\substack{0< k < m\\
(k,p)=1}}k,
\end{align*}
where the limit is for $m$ tending to $x$ $p$-adically in $\mathbb{Z}_{\ge 0}$.

We require several properties of $p$-adic Gamma function.
\begin{lem} (See  \cite[Section 11.6]{cohen-2007}.)
For any odd prime $p$ and $x,y\in \mathbb{Z}_p$, we have
\begin{align}
&\Gamma_p(1)=-1,\label{bb1}\\
&\Gamma_p(x)\Gamma_p(1-x)=(-1)^{s_p(x)},\label{bb2}\\
&\Gamma_p(x)\equiv \Gamma_p(y)\pmod{p}\quad\text{for $x\equiv y\pmod{p}$},\label{bb4}
\end{align}
where $s_p(x)\in \{1,2,\cdots,p\}$ with $s_p(x)\equiv x\pmod{p}$.
\end{lem}

\begin{lem} (See \cite[Lemma 17, (4)]{lr-2016}.)
Let $p$ be an odd prime. If $a\in \mathbb{Z}_p, n\in \mathbb{N}$
such that none of $a,a+1,\cdots,a+n-1$ in $p\mathbb{Z}_p$, then
\begin{align}
(a)_n=(-1)^n\frac{\Gamma_p(a+n)}{\Gamma_p(a)}.\label{bb5}
\end{align}
\end{lem}

\section{Proof of Theorem \ref{t1}}
Let $\omega$ be any primitive $3$th root of unity.
Letting $a=\frac{1}{2}, b=\frac{1-\omega p}{2}, c=\frac{1-\omega^2p}{2}, k=\frac{3}{2}, m=\frac{p-1}{2}$ in \cite[(1), page 32]{bailey-b-1964}, we obtain
\begin{align}
&\pFq{4}{3}{\frac{1}{2},&\frac{1-\omega p}{2},&\frac{1-\omega^2 p}{2},&\frac{1-p}{2}}{&1+\frac{\omega p}{2},&1+\frac{\omega^2 p}{2},&1+\frac{p}{2}}{1}\notag\\[10pt]
&=\frac{p\left(\frac{1}{2}\right)_\frac{p-1}{2}\left(\frac{1-p}{2}\right)_\frac{p-1}{2}}{\left(1+\frac{\omega p}{2}\right)_\frac{p-1}{2}\left(1+\frac{\omega^2 p}{2}\right)_\frac{p-1}{2}}
\pFq{4}{3}{\frac{1}{2},&\frac{1-\omega}{2},&\frac{1-\omega^2 p}{2},&\frac{1-p}{2}}{&1,&\frac{3}{4},&\frac{5}{4}}{1}.\label{b1}
\end{align}

By the fact that
\begin{align*}
(u+vp)(u+vp\omega)(u+vp\omega^2)=u^3+v^3p^3,
\end{align*}
we have
\begin{align}
(u+vp)_k(u+vp\omega)_k(u+vp\omega^2)_k\equiv (u)_k^3\pmod{p^3}.\label{b2}
\end{align}
It follows from \eqref{b1} and \eqref{b2} that
\begin{align}
&\pFq{4}{3}{\frac{1}{2},&\frac{1}{2},&\frac{1}{2},&\frac{1}{2}}{&1,&1,&1}{1}_{\frac{p-1}{2}}\notag\\[10pt]
&\equiv\frac{p\left(\frac{1}{2}\right)_\frac{p-1}{2}\left(\frac{1-p}{2}\right)_\frac{p-1}{2}}{\left(1+\frac{\omega p}{2}\right)_\frac{p-1}{2}\left(1+\frac{\omega^2 p}{2}\right)_\frac{p-1}{2}}
\pFq{4}{3}{\frac{1}{2},&\frac{1}{2},&\frac{1}{2},&\frac{1}{2}}{&1,&\frac{3}{4},&\frac{5}{4}}{1}_{\frac{p-1}{2}}\pmod{p^3}.
\label{b3}
\end{align}
In order to prove \eqref{a2}, by \eqref{a1} and \eqref{b3} it suffices to show that
\begin{align}
\frac{\left(\frac{1}{2}\right)_\frac{p-1}{2}\left(\frac{1-p}{2}\right)_\frac{p-1}{2}}{\left(1+\frac{\omega p}{2}\right)_\frac{p-1}{2}\left(1+\frac{\omega^2 p}{2}\right)_\frac{p-1}{2}}\equiv 1\pmod{p^3}.\label{b4}
\end{align}

By \eqref{b2}, we have
\begin{align*}
\left(1+\frac{p}{2}\right)_\frac{p-1}{2}\left(1+\frac{\omega p}{2}\right)_\frac{p-1}{2}\left(1+\frac{\omega^2 p}{2}\right)_\frac{p-1}{2}\equiv \left(1\right)^3_{\frac{p-1}{2}}\pmod{p^3},
\end{align*}
and so
\begin{align*}
\frac{\left(\frac{1}{2}\right)_\frac{p-1}{2}\left(\frac{1-p}{2}\right)_\frac{p-1}{2}}{\left(1+\frac{\omega p}{2}\right)_\frac{p-1}{2}\left(1+\frac{\omega^2 p}{2}\right)_\frac{p-1}{2}}
\equiv
\frac{\left(1+\frac{p}{2}\right)_\frac{p-1}{2}\left(\frac{1}{2}\right)_\frac{p-1}{2}\left(\frac{1-p}{2}\right)_\frac{p-1}{2}}
{\left(1\right)^3_{\frac{p-1}{2}}}\pmod{p^3}.
\end{align*}
Furthermore, we have
\begin{align*}
\left(\frac{1-p}{2}\right)_\frac{p-1}{2}=(-1)^{\frac{p-1}{2}}(1)_{\frac{p-1}{2}},
\end{align*}
and
\begin{align*}
\left(\frac{1}{2}\right)_\frac{p-1}{2}=(-1)^{\frac{p-1}{2}}\left(1-\frac{p}{2}\right)_{\frac{p-1}{2}}.
\end{align*}
Thus,
\begin{align}
\frac{\left(\frac{1}{2}\right)_\frac{p-1}{2}\left(\frac{1-p}{2}\right)_\frac{p-1}{2}}{\left(1+\frac{\omega p}{2}\right)_\frac{p-1}{2}\left(1+\frac{\omega^2 p}{2}\right)_\frac{p-1}{2}}\equiv
\frac{\left(1+\frac{p}{2}\right)_{\frac{p-1}{2}}\left(1-\frac{p}{2}\right)_{\frac{p-1}{2}}}
{\left(1\right)_{\frac{p-1}{2}}^2}\pmod{p^3}.\label{b5}
\end{align}

We next evaluate the product on the right-hand side of \eqref{b5} modulo $p^4$:
\begin{align*}
\frac{\left(1+\frac{p}{2}\right)_{\frac{p-1}{2}}\left(1-\frac{p}{2}\right)_{\frac{p-1}{2}}}
{\left(1\right)_{\frac{p-1}{2}}^2}=\prod_{j=1}^{\frac{p-1}{2}}\left(1-\frac{p^2}{4j^2}\right).
\end{align*}
From the following Taylor expansion:
\begin{align*}
\prod_{j=1}^{\frac{p-1}{2}}(a_j+b_jx^2)=\prod_{j=1}^{\frac{p-1}{2}}a_j\cdot\left(1+x^2\sum_{j=1}^{\frac{p-1}{2}}\frac{b_j}{a_j}\right)
+\mathcal{O}(x^4),
\end{align*}
we deduce that
\begin{align*}
\frac{\left(1+\frac{p}{2}\right)_{\frac{p-1}{2}}\left(1-\frac{p}{2}\right)_{\frac{p-1}{2}}}
{\left(1\right)_{\frac{p-1}{2}}^2}\equiv 1-\frac{p^2}{4}\sum_{j=1}^{\frac{p-1}{2}}\frac{1}{j^2}\pmod{p^4}.
\end{align*}
By Wolstenholme's theorem, we have
\begin{align*}
\sum_{j=1}^{\frac{p-1}{2}}\frac{1}{j^2}\equiv \frac{1}{2}\left(\sum_{j=1}^{\frac{p-1}{2}}\frac{1}{j^2}
+\sum_{j=1}^{\frac{p-1}{2}}\frac{1}{(p-j)^2}\right)=\frac{1}{2}\sum_{j=1}^{p-1}\frac{1}{j^2}\equiv 0\pmod{p}.
\end{align*}
It follows that
\begin{align}
\frac{\left(1+\frac{p}{2}\right)_{\frac{p-1}{2}}\left(1-\frac{p}{2}\right)_{\frac{p-1}{2}}}
{\left(1\right)_{\frac{p-1}{2}}^2}\equiv 1\pmod{p^3}.\label{b6}
\end{align}
Combining \eqref{b5} and \eqref{b6}, we complete the proof of \eqref{b4}.

\section{A human proof of Theorem \ref{t2}}
Letting $a=\frac{1}{2},x=2n+\frac{3}{2}$ in \cite[(14.1)]{Whipple-plms-1926}, we obtain
\begin{align}
&\pFq{6}{5}{\frac{5}{4},&\frac{1}{2},&-2n-1,&2n+2,&\frac{1}{2}+y,&\frac{1}{2}-y}
{&\frac{1}{4},&2n+\frac{5}{2},&-2n-\frac{1}{2},&1-y,&1+y}{-1}\notag\\[10pt]
&=\frac{\pi \Gamma\left(-2n-\frac{1}{2}\right)\Gamma\left(2n+\frac{5}{2}\right)
\Gamma\left(1+y\right)\Gamma\left(1-y\right)}
{\Gamma\left(\frac{1}{2}\right)\Gamma\left(\frac{3}{2}\right)
\Gamma\left(n+\frac{y}{2}+\frac{3}{2}\right)
\Gamma\left(n-\frac{y}{2}+\frac{3}{2}\right)
\Gamma\left(-n+\frac{y}{2}\right)
\Gamma\left(-n-\frac{y}{2}\right)}.\label{new-1}
\end{align}

Note that
\begin{align}
\Gamma\left(-2n-\frac{1}{2}\right)\Gamma\left(2n+\frac{5}{2}\right)
&\overset{\eqref{e1}}{=}\left(2n+\frac{3}{2}\right)\Gamma\left(-2n-\frac{1}{2}\right)
\Gamma\left(2n+\frac{3}{2}\right)\notag\\[10pt]
&\overset{\eqref{e3}}{=}\frac{(4n+3)\pi}{2\sin\left(\left(2n+\frac{3}{2}\right)\pi\right)}\notag\\[10pt]
&=-\frac{(4n+3)\pi}{2},\label{new-2}
\end{align}
and
\begin{align}
\Gamma\left(\frac{1}{2}\right)\Gamma\left(\frac{3}{2}\right)
=\frac{\pi}{2}.\label{new-3}
\end{align}
Also,
\begin{align}
\Gamma(1+y)\Gamma(1-y)\overset{\eqref{e4}}{=}\frac{1}{\pi}\Gamma\left(\frac{1+y}{2}\right)
\Gamma\left(\frac{1-y}{2}\right)\Gamma\left(1+\frac{y}{2}\right)
\Gamma\left(1-\frac{y}{2}\right).\label{new-4}
\end{align}
Substituting \eqref{new-2}--\eqref{new-4} into the right-hand side of \eqref{new-1} gives
\begin{align}
&\pFq{6}{5}{\frac{5}{4},&\frac{1}{2},&-2n-1,&2n+2,&\frac{1}{2}+y,&\frac{1}{2}-y}
{&\frac{1}{4},&2n+\frac{5}{2},&-2n-\frac{1}{2},&1-y,&1+y}{-1}\notag\\[10pt]
&=-(4n+3)\frac{\Gamma\left(\frac{1+y}{2}\right)
\Gamma\left(\frac{1-y}{2}\right)\Gamma\left(1+\frac{y}{2}\right)
\Gamma\left(1-\frac{y}{2}\right)}{\Gamma\left(n+\frac{y}{2}+\frac{3}{2}\right)
\Gamma\left(n-\frac{y}{2}+\frac{3}{2}\right)
\Gamma\left(-n+\frac{y}{2}\right)
\Gamma\left(-n-\frac{y}{2}\right)}\notag\\[10pt]
&\overset{\eqref{e2}}{=}-(4n+3)\frac{\left(-n+\frac{y}{2}\right)_{n+1}\left(-n-\frac{y}{2}\right)_{n+1}}
{\left(\frac{1+y}{2}\right)_{n+1}\left(\frac{1-y}{2}\right)_{n+1}}\notag\\[10pt]
&=-(4n+3)\frac{\left(\frac{y}{2}\right)_{n+1}\left(-n+\frac{y}{2}\right)_{n+1}}
{\left(-n+\frac{y-1}{2}\right)_{2n+2}}.\label{c1}
\end{align}

Let $i$ be any primitive $4$th root of unity.
Setting $n=\frac{p-3}{4}$ and $y=-\frac{ip}{2}$ in \eqref{c1} yields
\begin{align}
&\pFq{6}{5}{\frac{5}{4},&\frac{1}{2},&\frac{1-p}{2},&\frac{1+p}{2},&\frac{1-ip}{2},&\frac{1+ip}{2}}
{&\frac{1}{4},&1-\frac{p}{2},&1+\frac{p}{2},&1-\frac{ip}{2},&1+\frac{ip}{2}}{-1}\notag\\[10pt]
&=-\frac{p\left(-\frac{ip}{4}\right)_{\frac{p+1}{4}}
\left(\frac{3-(i+1)p}{4}\right)_{\frac{p+1}{4}}}
{\left(\frac{1-(i+1)p}{4}\right)_{\frac{p+1}{2}}}.
\label{c3}
\end{align}
By the fact that
\begin{align*}
(u+vp)(u-vp)(u+vpi)(u-vpi)=u^4-v^4p^4,
\end{align*}
we have
\begin{align}
(u+vp)_k(u-vp)_k(u+vpi)_k(u-vpi)_k\equiv (u)_k^4\pmod{p^4}.\label{c4}
\end{align}
In order to prove \eqref{a4}, by \eqref{c3} and \eqref{c4} it suffices to show that
\begin{align}
-\frac{p\left(-\frac{ip}{4}\right)_{\frac{p+1}{4}}
\left(\frac{3-(i+1)p}{4}\right)_{\frac{p+1}{4}}}
{\left(\frac{1-(i+1)p}{4}\right)_{\frac{p+1}{2}}}
\equiv -\frac{p^3}{16}\Gamma_p\left(\frac{1}{4}\right)^4\pmod{p^4}.\label{c5}
\end{align}

Note that
\begin{align*}
\left(-\frac{ip}{4}\right)_{\frac{p+1}{4}}
\left(\frac{3-(i+1)p}{4}\right)_{\frac{p+1}{4}}
=-\frac{p^2}{16}\prod_{j=1}^{\frac{p-3}{4}}\left(-\frac{p^2}{16}-j^2\right),
\end{align*}
and
\begin{align*}
\left(\frac{1-(i+1)p}{4}\right)_{\frac{p+1}{2}}=\prod_{j=1}^{\frac{p+1}{4}}
\left(-\frac{p^2}{16}-\left(-\frac{1}{2}+j\right)^2\right).
\end{align*}
It follows that
\begin{align*}
-\frac{p\left(-\frac{ip}{4}\right)_{\frac{p+1}{4}}
\left(\frac{3-(i+1)p}{4}\right)_{\frac{p+1}{4}}}
{\left(\frac{1-(i+1)p}{4}\right)_{\frac{p+1}{2}}}
&=\frac{p^3}{16}\cdot\frac{\prod_{j=1}^{\frac{p-3}{4}}\left(-\frac{p^2}{16}-j^2\right)}
{\prod_{j=1}^{\frac{p+1}{4}}
\left(-\frac{p^2}{16}-\left(-\frac{1}{2}+j\right)^2\right)}\\
&\equiv -\frac{p^3}{16}\cdot\frac{\left(1\right)_{\frac{p-3}{4}}^2}
{\left(\frac{1}{2}\right)_{\frac{p+1}{4}}^2}\pmod{p^4}\\
&\overset{\eqref{bb5}}{=} -\frac{p^3}{16}\cdot
\frac{\Gamma_p\left(\frac{p+1}{4}\right)^2\Gamma_p\left(\frac{1}{2}\right)^2}
{\Gamma_p\left(1\right)^2\Gamma_p\left(\frac{p+3}{4}\right)^2}.
\end{align*}
Furthermore, by \eqref{bb1}--\eqref{bb4}, we have
\begin{align*}
-\frac{p\left(-\frac{ip}{4}\right)_{\frac{p+1}{4}}
\left(\frac{3-(i+1)p}{4}\right)_{\frac{p+1}{4}}}
{\left(\frac{1-(i+1)p}{4}\right)_{\frac{p+1}{2}}}
&\equiv -\frac{p^3}{16}\cdot
\frac{\Gamma_p\left(\frac{1}{4}\right)^2}
{\Gamma_p\left(\frac{3}{4}\right)^2}\pmod{p^4}\\
&=-\frac{p^3}{16}\Gamma_p\left(\frac{1}{4}\right)^4.
\end{align*}
This completes the proof of \eqref{c5}.

\vskip 5mm \noindent{\bf Acknowledgments.} The author would like to thank Dr. Chen Wang for his helpful comments on this paper.
This work was supported by the National Natural Science Foundation of China (grant 11801417).

\end{document}